\documentclass[USletter, 10pt, conference]{ieeeconf}

\IEEEoverridecommandlockouts
\usepackage{cite}
\usepackage[acronym]{glossaries}
\usepackage{amsmath,amssymb,amsfonts}
\usepackage{algorithmic}
\usepackage{graphicx}
\usepackage{xcolor}
\def\BibTeX{{\rm B\kern-.05em{\sc i\kern-.025em b}\kern-.08em
    T\kern-.1667em\lower.7ex\hbox{E}\kern-.125emX}}

\usepackage{units}
\usepackage{mathtools} 
\usepackage{amsthm}
\usepackage{lettrine}

\usepackage{subcaption}
\usepackage{scrextend} 

\usepackage[english]{babel}
\usepackage{cite}

\usepackage{array}
\usepackage{color}
\usepackage{amssymb}
\usepackage{multicol}

\DeclareUnicodeCharacter{2212}{-}

\usepackage{todonotes}

\usepackage{tikz}
\usetikzlibrary{arrows,calc}
\usepackage{pgfplots}
\usepgfplotslibrary{groupplots}

\usepackage{amsthm}
\usepackage{epstopdf}

\newcounter{thm}
\newtheoremstyle{mystyle}
  {}
  {}
  {}
  {}
  {\bfseries\color{black}}
  {.}
  {\newline}
  {\thmname{#1}\thmnumber{ #2: }\normalfont\color{black}\thmnote{ \textit{(#3)}}}
\theoremstyle{mystyle}
\newtheorem{prob}[thm]{Problem}

\newtheorem{theorem}{Theorem}[section]

\usepackage{subfiles}

\usepackage[ruled]{algorithm2e}
\usepackage[usestackEOL]{stackengine}
\usepackage{makecell}
\usepackage{tabularray}

\usepackage{adjustbox}
\usepackage{nccmath}

\usepackage{url}
\newacronym{acr:cvt}{CVT}{continuously variable transmission}
\newacronym{acr:CoG}{CoG}{center of gravity}
\newacronym{acr:CoP}{CoP}{center of pressure}
\newacronym{acr:CP}{CP}{connection point}
\newacronym{acr:dp}{DP}{dynamic programming}
\newacronym{acr:ecms}{ECMS}{equivalent consumption minimization strategies}
\newacronym{acr:eltms}{ELTMS}{equivalent lap time minimization strategies}
\newacronym{acr:em}{EM}{electric motor}
\newacronym{acr:es2k}{ES2K}{Energy Storage to Kinetic}
\newacronym{acr:F1}{F1}{Formula 1}
\newacronym{acr:FIA}{FIA}{F\'{e}d\'{e}ration Internationale de l'Automobile}
\newacronym{acr:fgt}{FGT}{fixed-gear transmission}
\newacronym{acr:FD}{FD}{final drive}
\newacronym{acr:ice}{ICE}{internal combustion engine}
\newacronym{acr:k2es}{K2ES}{Kinetic to Energy Storage}
\newacronym{acr:mgu}{MGU}{motor generator unit}
\newacronym{acr:mguh}{MGU-H}{motor generator unit heat}
\newacronym{acr:mguk}{MGU-K}{motor generator unit kinetic}
\newacronym{acr:mpc}{MPC}{model predictive control}
\newacronym{acr:MM}{MM}{mounted motor}
\newacronym{acr:HM}{HM}{hub motor}

\newacronym[description={Energy Management Strategy}, \glslongpluralkey={Energy Management Strategies},\glsshortpluralkey={EMSs}]{EMS}{EMS}{Energy Management Strategy}%
\newacronym{acr:pmp}{PMP}{Pontryagin's Minimum Principle}
\newacronym{acr:pu}{PU}{power unit}
\newacronym[description={Powertrain Operation}, \glslongpluralkey={Powertrain Operations},\glsshortpluralkey={POs}]{acr:PO}{PO}{Powertrain Operation}%
\newacronym{acr:rmse}{RMSE}{root-mean-square error}
\newacronym{acr:socp}{SOCP}{second-order cone program}
\newacronym{acr:soe}{SoE}{State of Energy}

\newcommand{\pushright}[1]{\ifmeasuring@#1\else\omit\hfill$\displaystyle#1$\fi\ignorespaces}
\newcommand{\pushleft}[1]{\ifmeasuring@#1\else\omit$\displaystyle#1$\hfill\fi\ignorespaces}
\makeatother
\newif\ifmargincomments 
\margincommentstrue

\newif\ifextendedversion 
\extendedversionfalse

\ifmargincomments

\else

\fi
\usetikzlibrary{external}
\begin{document}

%
\title{ \LARGE \bf Two-dimensional Spatial Optimization for Electric Motorcycle Powertrain Elements using Mixed-integer Programming}

%
%
%

\author{Jorn van Kampen$^\star$, Chun-Cheng Huang$^\star$, Mauro Salazar
	\thanks{$^\star$ These authors contributed equally to this work.}\thanks{ Jorn van Kampen, Chun-Cheng Huang and Mauro Salazar are with the Control Systems Technology section, Department of Mechanical Engineering, Eindhoven University of Technology (TU/e), Eindhoven, 5600 MB, The Netherlands.
		E-mails: {\tt\footnotesize j.h.e.v.kampen@tue.nl}, {\tt\footnotesize chun.cheng0702@gmail.com}, {\tt\footnotesize m.r.u.salazar@tue.nl} 
}}

\maketitle

\begin{abstract}
This study presents a framework for optimizing the two-dimensional (2D) placement of electric motorcycle powertrain elements, accounting for the position, the orientation and geometric irregularities. 
Specifically, we construct a 2D placement model at the component level in which we include near-continuous rotation of components and allow for irregular subsystem geometries to make optimal use of the limited design space.
Second, we introduce linearization techniques for the trigonometric constraints and formulate the placement problem as a mixed-integer quadratic program (MIQP). 
Finally, we demonstrate our framework on two electric motorcycle powertrain topologies and study the influence of the geometry complexity on the placement solutions. 
The results show that gradually increasing complexity leads to more manageable computation times and higher the complexity solution improves handling performance by 2.5\% compared to the benchmark placement found in existing electric motorcycles.

\end{abstract}


\section{Introduction}
\begin{figure}[t]
	\centering
	\includegraphics[trim=2.2cm 0cm 0cm 0cm,clip,width=1\columnwidth]{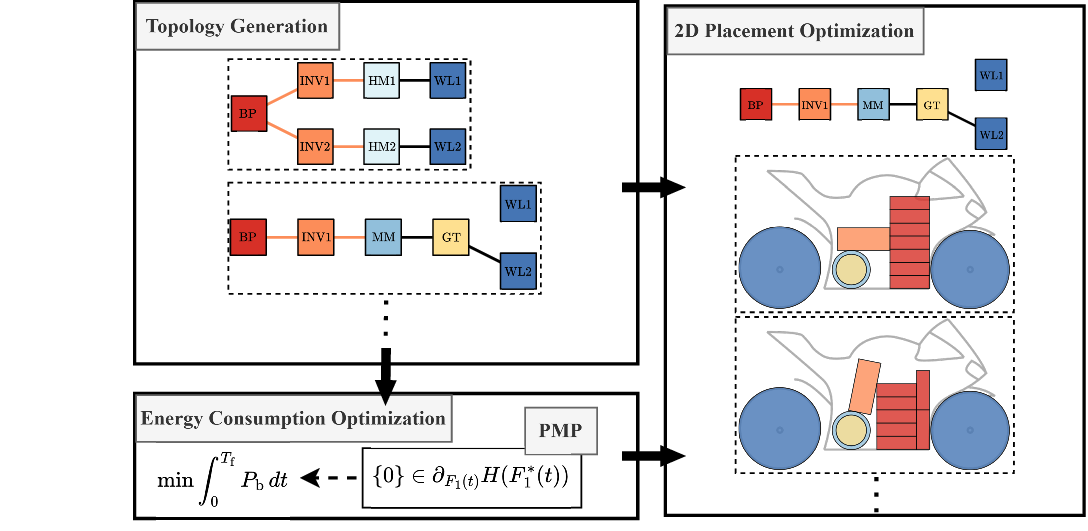}
	\caption{Overview of the integrated framework for electric motorcycles. The workflow begins with the generation of a set of feasible topologies. In the second stage, the optimal energy consumption of each topology is determined, incorporating a modified Pontryagin’s Minimum Principle to derive the optimal control policy. The final stage optimizes the component placement to enhance handling, while maintaining optimal energy consumption.}
	\label{fig: Fig1_demo}
\end{figure}

\lettrine[findent=2pt]{\textbf{E}}{lectric} vehicles represent a significant step towards sustainable transportation, offering a cleaner alternative to conventional fuel-powered vehicles. While electric cars have received considerable attention and development, they are still relatively expensive and inefficient in terms of road space. To overcome these challenges, electric motorcycles could form an alternative, since they are considerably cheaper and more compact. 
As with electric cars, the design and operation of the powertrain has a significant influence on the Total Cost of oownership (TCO). Consequently, the powertrain architecture, its component specifications and the energy management strategy should be carefully optimized, motivating the need for integrated approaches~\cite{WijknietHofman2018}. A common example of such an approach is shown on the left half in Fig.\ref{fig: Fig1_demo}, whereby the powertrain topology and energy consumption are optimized.  
However, electric motorcycles also come with their own unique set of challenges.
Their compact and lightweight design complicates component placement due to limited space and makes the motorcycle's \gls{acr:CoG} highly sensitive to the positioning of components. The \gls{acr:CoG} is a critical factor influencing static weight distribution, dynamic weight transfer and overall handling. Suboptimal \gls{acr:CoG} positioning can compromise dynamic performance, increase energy consumption, and degrade riding experience. 
Furthermore, the limited design space requires the exploration of irregular geometries for the powertrain elements to maximize the spatial utilization.
This highlights the need for a design strategy capable of computing the placement for arbitrary powertrain topologies with irregular geometries. 
Against this backdrop, this paper presents an optimization framework to compute the optimal powertrain element placement for electric motorcycles with arbitrary powertrain architectures. The proposed method builds upon existing integrated design approaches, as shown on the right side of Fig.~\ref{fig: Fig1_demo}.


\paragraph*{Related Literature}




The design optimization of vehicle powertrain systems often starts from the topology selection~\cite{WijknietHofman2018, ChenouardHartmannEtAl2016, ZhaoTangEtAl2022}. 
Several works have extended these topology selection algorithms with spatial optimization~\cite{KabalanVinotEtAl2020, BerxGadeyneEtAl2014, PiacentiniCheongEtAl2020}. 
However, these studies often focus on smaller subsystems within vehicle powertrains, such as gearboxes, and typically constraint certain degrees of freedom in the design, such as component orientation or relative positioning.
While we addressed this in our previous work~\cite{KampenSalazarEtAl2023,KampenSalazarEtAl2025}, we restricted element geometries to cuboidal shapes or composite polycubes constructed from combined cuboid primitives.
Other studies incorporated arbitrary element geometries by directly using CAD models~\cite{HofstetterLechleitnerEtAl2018,FilomenoAhmadEtAl2021}, yet they applied gradient-free optimization methods, which provide no information on the solution quality for non-convex problems. 

For more general two- or three-dimensional placement problems, the floorplanning problem --- where a set of rectangles have to be arranged such that their bounding box is minimized --- serves as a basis for many layout design problems~\cite{Laignel2021, Tresca2023}. However, heterogeneity in element geometry beyond variety in dimensions of the rectangles is often disregarded.
Finally, arbitrary element shapes in general placement problems are considered in~\cite{CaganShimadaEtAl2002, BelloEtAl2024, PeddadaEtAl2022}, but again there are no measures for the quality of the solution, and the obtained solutions are sensitive to the initial guess. 

In conclusion, to the best of the authors' knowledge, there are currently no studies that optimize the spatial powertrain design for irregular geometries and non-orthogonal orientations whilst providing a qualitative measure of the optimality of the solution.  

%
%

\paragraph*{Statement of Contributions}
This paper presents an optimization framework for the two-dimensional design of electric motorcycle powertrains to maximize handling performance for arbitrary topologies, ultimately aiming to build upon existing integrated design frameworks. 
We base ourselves on our previous works regarding spatial optimization~\cite{KampenSalazarEtAl2023, KampenSalazarEtAl2025} and extend them to allow for near-continuous rotation of components and irregular geometries of subsystems, where we introduce clusters to control the design complexity of the subsystems. 
Subsequently, we introduce various linearization techniques and frame the problem as a Mixed-integer Quadratic Program (MIQP) that can be solved by off-the-shelf algorithms and provides a qualitative measure of the solutions quality through the optimality gap. 
Finally, we showcase our framework for two powertrain architectures and study the effect of the amount of allowed clusters on the resulting powertrain design.
%

\paragraph*{Organization}
The remainder of this paper is structured as follows: Section~\ref{sec:method_2D} introduces the spatial placement problem, including the extensions to non-orthogonal orientations and irregular geometries.
In Section~\ref{sec:results}, we demonstrate our framework for two electric motorcycle topologies generated by the existing integrated framework. 
Finally, we draw the conclusion in Section and provide an outlook in~\ref{sec: conclusion}.




\section{Two-dimensional element placement}
\label{sec:method_2D}
This section derives the two-dimensional optimal element placement for a given motorcycle powertrain topology $\mathbf{T}$. The goal is to find the optimal position and orientation of all elements, and optimal configuration of subsystems, inspired by our previous work on two-dimensional~\cite{KampenSalazarEtAl2023} and three-dimensional~\cite{KampenSalazarEtAl2025} spatial optimization. 
In contrast to our previous work, we now include finer discretized orientations, introduce circular component shapes to the model and allow for more design freedom in the subsystem shape.
The current framework includes six element types for the motorcycle powertrain: mounted-motors (MM), hub-motors (HM), inverters (INV), a battery pack (BP), gearbox and transmission (GT), and wheels (WL). Each type has a maximum number of instances $N_\tau$, with each element denoted as $\tau_j$, $j\in\{1,\dots,N_\tau\}$, and its existence indicated by $\mathbf{n}_{\tau_j}\in\{0,1\}$. Some elements, such as the battery pack, are regarded as subsystems that contain internal elements, such as battery modules. In this work, the subsystem design is determined solely by the placement and orientation of the internal elements and is no longer restricted to be a bounding box.  

\subsection{Design Space and Shape Definition}
The design space represents the region designated for element placement and configuration, and is defined in a Cartesian coordinate system with the origin located at the rear wheel contact point. The longitudinal boundary is specified as $[\underline{x}, \overline{x}]$, and the vertical boundary is defined as $[\underline{y}, \overline{y}]$. In this framework, we consider two shapes for the components: rectangular for the inverters and battery modules, and circular for electric motors and the transmission. Rectangles have attributes: width $w$, height $h$, mass $m$, center point $(x,y)$, and angle $\theta$. Circles have attributes: radius $r$, mass $m$, and center point $(x,y)$, and projected angle $\theta_p$. The center point of each shape is placed within the design space $x\in[\underline{x},\overline{x}]$, $y\in[\underline{y},\overline{y}]$. For a rectangle, we introduce two unit vectors, $\vec{u}_{1} = (\mathrm{cos}(\theta), \mathrm{sin}(\theta))$ and $\vec{u}_{2} = (-\mathrm{sin}(\theta), \mathrm{cos}(\theta))$, and define the discretized orientation as
\begin{equation}
	\theta_k = 180^\circ \cdot \frac{k}{N_\mathrm{a}}, \: k \in \mathcal{K}_\mathrm{a}= \{0,\ldots,N_\mathrm{a}-1\},
	\label{eq: 2D-d-angle}
\end{equation} 
where $N_a$ (an even integer) represents the total number of discretized angles. 
For a circle, we consider a projected vector $\vec{u}_{\mathrm{p}} = (\mathrm{cos}(\theta_{\mathrm{p}}), \mathrm{sin}(\theta_{\mathrm{p}}))$, and define the discretized projected angle as 
\begin{equation}
	\theta_{k_\mathrm{p}} = 180^\circ \cdot \frac{k_\mathrm{p}}{N_\mathrm{pa}-1}, \: {k_\mathrm{p}} \in \mathcal{K}_\mathrm{pa}= \{0,\ldots,N_\mathrm{pa}-1\},
\end{equation}    
where $N_\mathrm{pa}$ is the number of projected angles.

For a rectangular object $\mathfrak{d}$, we define the angle decision set as $\mathcal{D}_{\mathfrak{d}} \coloneqq \{\delta_{k,\mathfrak{d}} \in \{0,1\} \mid \forall k \in \mathcal{K}_\mathrm{a}\}$, where $\delta_{k,\mathfrak{d}}$ is the angle decision variable. The object's angle is given by
\begin{equation}
	\theta_{\mathfrak{d}} = \sum_{k=0}^{N_\mathrm{a}-1}\delta_{k,\mathfrak{d}} \cdot \theta_k, \text{ subject\ to }\sum\mathcal{D}_{\mathfrak{d}}=1,
\end{equation} 
ensuring a single angle is selected. Similarly, for a circular object $\mathfrak{z}$, we define the projected angle decision set as $\mathcal{D}_{\mathrm{p},\mathfrak{z}} \coloneqq \{\delta_{k_\mathrm{p},\mathfrak{z}} \in \{0,1\} \mid \forall k_\mathrm{p} \in \mathcal{K}_\mathrm{pa}\}$, and the object's projected angle is
\begin{equation}
	\theta_{\mathrm{p},\mathfrak{z}} = \sum_{k=0}^{N_\mathrm{pa}-1}\delta_{k_\mathrm{p},\mathfrak{z}} \cdot \theta_{k_\mathrm{p}}, \text{ subject\ to }\sum\mathcal{D}_{\mathrm{p},\mathfrak{z}}=1.
\end{equation}

In the following content, we denote the attributes of an object's shape with a subscript, such as the width $w_\mathfrak{d}$.

\subsection{Overlap Prevention}
\begin{figure}[t]
	\begin{center} 
		\begin{subfigure}[b]{0.49\linewidth}
			\includegraphics[trim=1.8cm 0cm 0cm 0cm,clip,width=\textwidth]{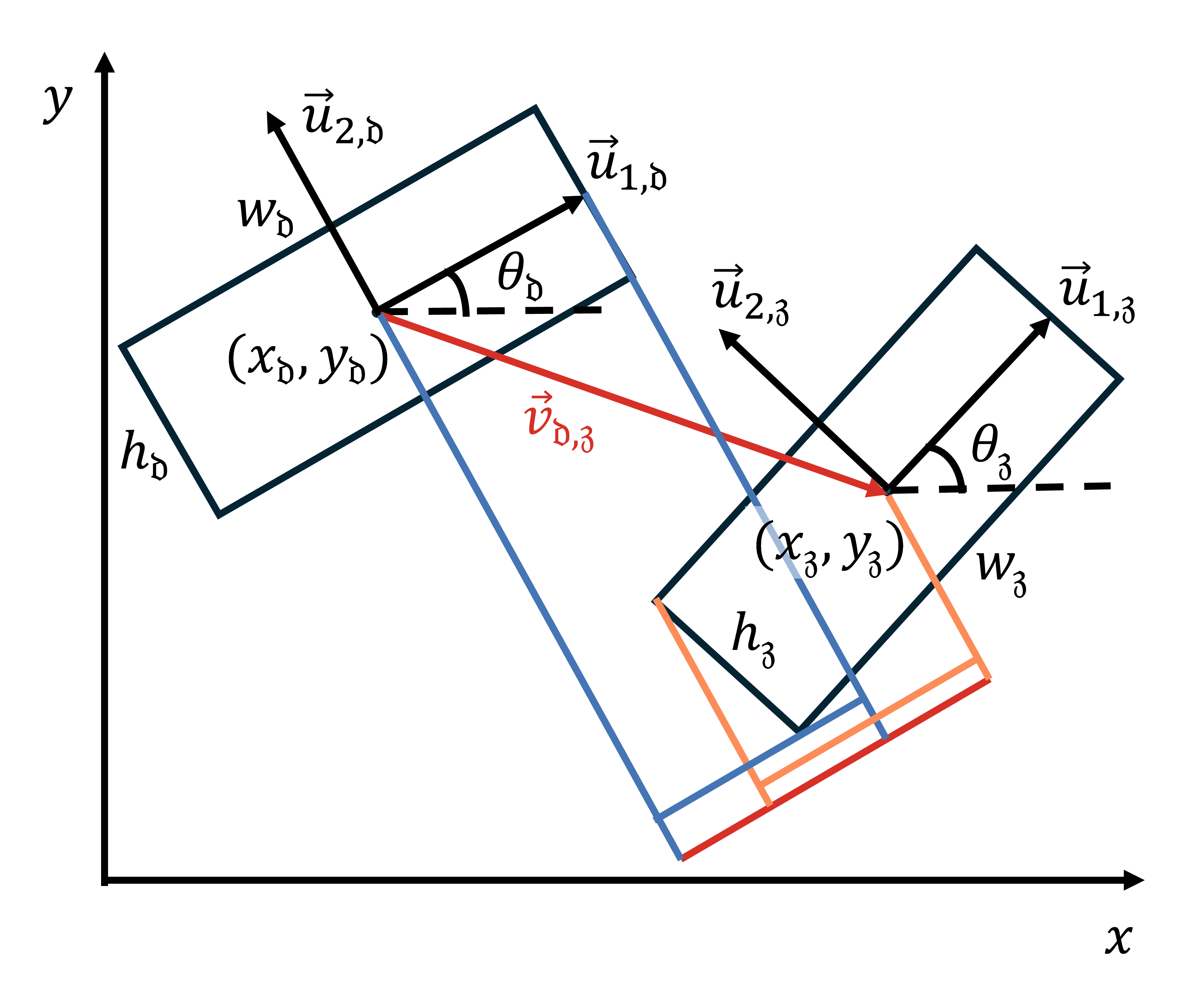}
			\caption{}
			\label{fig: SAT-demo-failed}
		\end{subfigure}
		\hfill
		\begin{subfigure}[b]{0.49\linewidth}
			\includegraphics[trim=1.8cm 0cm 0cm 0cm,clip,width=\textwidth]{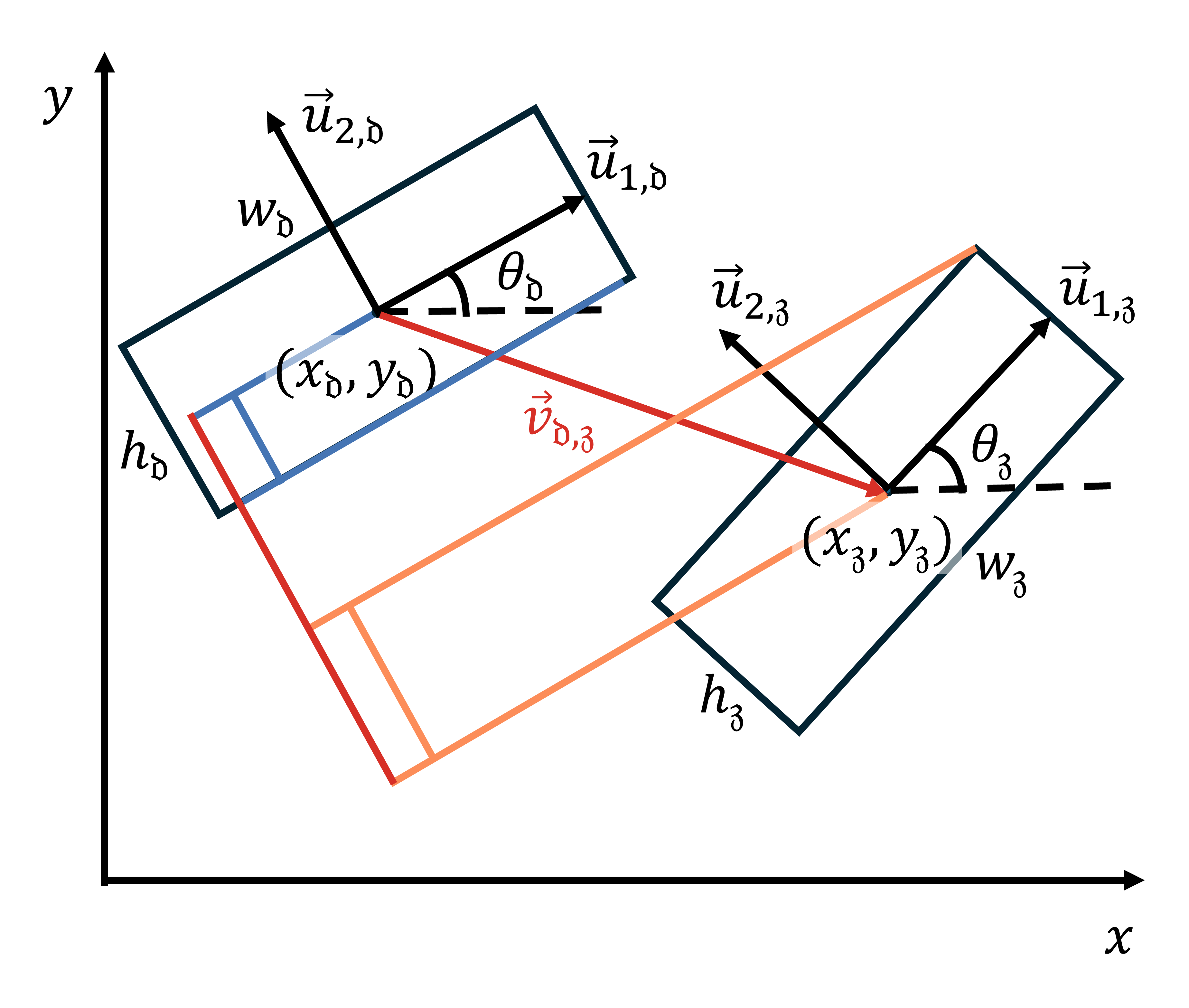}
			\caption{}
			\label{fig: SAT-demo-pass}
		\end{subfigure}
		\caption{Visualization of the SAT applied to two rectangular objects. In case (\subref{fig: SAT-demo-failed}), condition 1 in~\eqref{eq: 2D-SAT-rec} fails, as the length of the projected center vector $\vec{v}_{\mathfrak{d},\mathfrak{z}}$ onto the unit vector $\vec{u}_{1,\mathfrak{d}}$ is shorter than the projected lengths of both objects. In contrast, case (\subref{fig: SAT-demo-pass}) illustrates that condition 2 is satisfied, meaning that the objects do not overlap.} 
		\label{fig: SAT-demo}
	\end{center}
\end{figure}

For the resulting element placement to be feasible, we require that none of the elements overlap with one other. 
%
We do this by applying the Separate Axis Theorem (SAT)~\cite{Liang2015} to every pair of objects $\mathfrak{d}, \mathfrak{z} \in \mathbf{T}$. 
For two rectangular objects, the shapes do not overlap if either of the conditions is satisfied:
\begin{equation}
	\begin{cases}
		1. \:|\vec{v}_{\mathfrak{d},\mathfrak{z}}\cdot \vec{u}_{1,\mathfrak{d}}| \geq \frac{w_{\mathfrak{d}}}{2} + |\frac{w_{\mathfrak{z}}}{2} \cdot \vec{u}_{1,\mathfrak{z}} \cdot \vec{u}_{1,\mathfrak{d}}| + |\frac{h_{\mathfrak{z}}}{2} \cdot \vec{u}_{2,\mathfrak{z}} \cdot \vec{u}_{1,\mathfrak{d}}| \\
		2. \:|\vec{v}_{\mathfrak{d},\mathfrak{z}}\cdot \vec{u}_{2,\mathfrak{d}}| \geq \frac{h_{\mathfrak{d}}}{2} + |\frac{w_{\mathfrak{z}}}{2} \cdot \vec{u}_{1,\mathfrak{z}} \cdot \vec{u}_{2,\mathfrak{d}}| + |\frac{h_{\mathfrak{z}}}{2} \cdot \vec{u}_{2,\mathfrak{z}} \cdot \vec{u}_{2,\mathfrak{d}}| \\
		3. \:|\vec{v}_{\mathfrak{d},\mathfrak{z}}\cdot \vec{u}_{1,\mathfrak{z}}| \geq \frac{w_{\mathfrak{z}}}{2} + |\frac{w_{\mathfrak{d}}}{2} \cdot \vec{u}_{1,\mathfrak{d}} \cdot \vec{u}_{1,\mathfrak{z}}| + |\frac{h_{\mathfrak{d}}}{2} \cdot \vec{u}_{2,\mathfrak{d}} \cdot \vec{u}_{1,\mathfrak{z}}| \\
		4. \:|\vec{v}_{\mathfrak{d},\mathfrak{z}}\cdot \vec{u}_{2,\mathfrak{z}}| \geq \frac{h_{\mathfrak{z}}}{2} + |\frac{w_{\mathfrak{d}}}{2} \cdot \vec{u}_{1,\mathfrak{d}} \cdot \vec{u}_{2,\mathfrak{z}}| + |\frac{h_{\mathfrak{d}}}{2} \cdot \vec{u}_{2,\mathfrak{d}} \cdot \vec{u}_{2,\mathfrak{z}}|,
	\end{cases}
	\label{eq: 2D-SAT-rec}
\end{equation}
where $\vec{v}_{\mathfrak{d},\mathfrak{z}} = (x_{\mathfrak{z}}-x_{\mathfrak{d}}, \: y_{\mathfrak{z}}-y_{\mathfrak{d}})$, is the center vector pointing from the center of object $\mathfrak{d}$ to the center of object $\mathfrak{z}$. The visualization of the SAT applied to two rectangular objects is shown in Fig.~\ref{fig: SAT-demo}. We can apply similar reasoning for a pair $\mathfrak{d}, \mathfrak{z}$ where $\mathfrak{d}$ is rectangular and $\mathfrak{z}$ is circular:
\begin{equation}
	\begin{aligned}
		&\vec{v}_{\mathfrak{d},\mathfrak{z}}^\mathrm{r} = 
		\begin{bmatrix}
			\mathrm{cos}(\theta_{\mathrm{r},\mathfrak{d}}) & \mathrm{sin}(\theta_{\mathrm{r},\mathfrak{d}})\\
			-\mathrm{sin}(\theta_{\mathrm{r},\mathfrak{d}}) & \mathrm{cos}(\theta_{\mathrm{r},\mathfrak{d}})
		\end{bmatrix} \cdot
		\begin{bmatrix}
			x_{\mathfrak{z}}-x_{\mathfrak{d}}\\
			y_{\mathfrak{z}}-y_{\mathfrak{d}}
		\end{bmatrix},\\
		&|\vec{v}_{\mathfrak{d},\mathfrak{z}}^\mathrm{r}\cdot\vec{u}_{\mathrm{p},\mathfrak{z}}| \geq r_{\mathfrak{z}} + \frac{w_{\mathfrak{d}}}{2}\cdot|\mathrm{cos}(\theta_{\mathrm{p},\mathfrak{z}})| + \frac{h_{\mathfrak{d}}}{2}\cdot|\mathrm{sin}(\theta_{\mathrm{p},\mathfrak{z}})|,
	\end{aligned}
	\label{eq: 2D-SAT-circle}
\end{equation} 
where $\vec{v}_{\mathfrak{d},\mathfrak{z}}^\mathrm{r}$ is the rotated center vector.

\subsection{Component and Cluster Formulation}
Elements in the powertrain system can be standalone components or composite subsystems. In the case of subsystems, we want to generate arbitrary shapes instead of rectangles to make optimal use of the limited design space present in motorcycles. 
To this end, we consider rectangular clusters of sub-modules $\tau_{j,\text{sub}}$ that should together represent the subsystems as a continuous shape.
The number of sub-modules in a subsystem is $N_{\text{sub},\tau_j}\ge 1$, while components have $N_{\text{sub},\tau_j}=0$. 

The rectangular clusters are captured by a set of all possible sub-module arrangements $\mathcal{S}_{\tau_j} \coloneq \left\{\tau_{j,\mathfrak{s}} | \mathfrak{s} \in \mathcal{I}_{\tau_j}\right\}$, where $\tau_{j,\mathfrak{s}}$ represents the cluster arrangement with identifier $\mathfrak{s}$. A rectangular arrangement consists of $n_{\mathrm{b},\tau_{j,\mathfrak{s}}}$ elements distributed over $n_{\mathrm{w},\tau_{j,\mathfrak{s}}}$ columns and $n_{\mathrm{h},\tau_{j,\mathfrak{s}}}$ rows such that
\begin{equation}
	\begin{aligned}
		 w_{\tau_{j,\mathfrak{s}}} &= n_{\mathrm{w},\tau_{j,\mathfrak{s}}}\cdot w_{\tau_{j,\mathrm{sub}}}, \quad 
		h_{\tau_{j,\mathfrak{s}}} = n_{\mathrm{h},\tau_{j,\mathfrak{s}}}\cdot h_{\tau_{j,\mathrm{sub}}}, \\
		 n_{\mathrm{b},\tau_{j,\mathfrak{s}}} &= n_{\mathrm{w},\tau_{j,\mathfrak{s}}} \cdot n_{\mathrm{h},\tau_{j,\mathfrak{s}}}, \quad n_{\mathrm{w},\tau_{j,\mathfrak{s}}}, n_{\mathrm{h},\tau_{j,\mathfrak{s}}} \in \mathbb{N}^+, \\
		 m_{\tau_{j,\mathfrak{s}}} &= n_{\mathrm{b},\tau_{j,\mathfrak{s}}}\cdot m_{\tau_{j,\mathrm{sub}}}\leq m_{\tau_j}, \\ 
			w_{\tau_{j,\mathfrak{s}}} &\leq (\overline{x}-\underline{x}), \quad  h_{\tau_{j,\mathfrak{s}}} \leq (\overline{y}-\underline{y}).\\
	\end{aligned}
\end{equation}

The component's shape is constructed using $N_{\mathrm{com},\tau_{j}}$ pre-arranged clusters, where $N_{\mathrm{com},\tau_{j}} \in \{1,\ldots,N_{\mathrm{sub},\tau_j}\}$ is the shape complexity of the component. To select clusters, we define the decision set as $\mathcal{D}_{\tau_{j,\mathfrak{c}_i}} \coloneqq \{\delta_{\tau_{j,\mathfrak{s},\mathfrak{c}_i}} \in \{0,1\} \mid \forall \tau_{j,\mathfrak{s}} \in \mathcal{S}_{\tau_j}\}$, where the cluster decision variable $\delta_{\tau_{j,\mathfrak{s},\mathfrak{c}_i}}$ ensures each cluster $\tau_{j,\mathfrak{c}_i}$ is selected from $\mathcal{S}_{\tau_j}$, and subject to
\begin{equation}
	\sum{\mathcal{D}_{\tau_{j,\mathfrak{c}_i}}} = 1,
\end{equation}
ensuring a single pre-arranged cluster is selected. We then calculate the attributes of each selected cluster  $\tau_{j,\mathfrak{c}_i}$ as a linear combination of the attributes of the pre-arranged clusters. For example,
\begin{equation}
	w_{\tau_{j,\mathfrak{c}_i}} = \sum_{\tau_{j,\mathfrak{s}}\in \mathcal{S}_{\tau_j}}w_{\tau_{j,\mathfrak{s}}}\cdot\delta_{\tau_{j,\mathfrak{s},\mathfrak{c}_i}}.
\end{equation}
Then, we define the set of clusters that form the shape of component $\tau_j$ as
\begin{equation}
	\mathcal{C}_{\tau_j} \coloneqq \{\tau_{j,\mathfrak{c}_i} \mid \forall i \in \{1,\ldots,N_{\mathrm{com},\tau_j}\}, N_{\mathrm{sub},\tau_j}\geq1\}.
\end{equation}
%

The last step is ensuring that the resulting subsystems are well-defined.
First, we ensure the mass equivalence through
\begin{equation}
	\sum_{i=1}^{N_{\mathrm{com},\tau_j}}m_{\tau_{j,\mathfrak{c}_i}} = m_{\tau_j}, \quad \tau_{j,\mathfrak{c}_i} \in \mathcal{C}_{\tau_j},    
\end{equation}
which ensures that the total mass of the clusters equals the mass of the subsystem. 
Second, to reduce complexity, we enforce every pair of clusters $\mathfrak{d},\mathfrak{z} \in \mathcal{C}_{\tau_j}$ within the subsystem to have the same orientation or $90^{\circ}$ difference as
\begin{equation}
	\begin{aligned}
		& \delta_{k, \mathfrak{d}} \leq \delta_{k, \mathfrak{z}} + \delta_{k_{\perp}, \mathfrak{z}}, \\
		& k_{\perp} = \mathrm{mod}(k + \frac{N_\mathrm{a}}{2},N_\mathrm{a}), \quad k,k_{\perp}\in \mathcal{K}_\mathrm{a}, \\
	\end{aligned}
\end{equation}
where $k_{\perp}$ is the index of the angle that differs $90^\circ$ from the angle with index $k$, and $\mathrm{mod}(a,b)$ represents the remainder when integer $a$ is divided by integer $b$. 
To analyze the exact orientation relationship, we introduce the orientation indicator $\delta_{\mathrm{o},\mathfrak{d},\mathfrak{z}} \in \{0,1\}$ for each pair of clusters, defined as 
\begin{equation}
	(1-\delta_{\mathrm{o},\mathfrak{d},\mathfrak{z}}) \leq
	\sum_{k=0}^{N_\mathrm{a}-1}(\delta_{k, \mathfrak{d}}-\delta_{k, \mathfrak{z}})^2 \leq 2\cdot (1-\delta_{\mathrm{o},\mathfrak{d},\mathfrak{z}}).
\end{equation}
The indicator $\delta_{\mathrm{o},\mathfrak{d},\mathfrak{z}} = 1$ represents the two clusters have the identical orientation, while $\delta_{\mathrm{o},\mathfrak{d},\mathfrak{z}} = 0$ indicates their orientations differ by $90^\circ$.

Finally, to ensure that selected clusters form a continuous geometry, we propose relative position constraints with a big-M formulation~\cite{RichardsHow2005} as follows:
\begin{equation}
	\begin{aligned}
		& \vec{v}_{\mathfrak{d},\mathfrak{z}}\cdot \vec{u}_{1,\mathfrak{d}} \leq \left|\frac{w_{\mathfrak{d}} + l_{1,\mathfrak{z}}}{2} - \epsilon \cdot (1-\delta_{\mathrm{s}, \mathfrak{d},\mathfrak{z}}) + M \cdot \Delta_\mathrm{o}\right|, \\
		& \vec{v}_{\mathfrak{d},\mathfrak{z}}\cdot \vec{u}_{2,\mathfrak{d}} \leq \left|\frac{h_{\mathfrak{d}} + l_{2,\mathfrak{z}}}{2} - \epsilon \cdot \delta_{\mathrm{s}, \mathfrak{d},\mathfrak{z}} + M \cdot \Delta_\mathrm{o}\right|, \\
		& \delta_{\mathrm{s}, \mathfrak{d},\mathfrak{z}} \in \{0,1\}, \quad \forall \Delta_\mathrm{o} \in \{1-\delta_{\mathrm{o},\mathfrak{d},\mathfrak{z}}, \: \delta_{\mathrm{o},\mathfrak{d},\mathfrak{z}}\},
	\end{aligned}
	\label{eq: 2D-relative-placement}
\end{equation}
where $\delta_{\mathrm{s}, \mathfrak{d},\mathfrak{z}}$ is the relative position indicator, $\epsilon$ is the tolerance, and $M$ is a large constant number. The lengths $l_{1,\mathfrak{z}}$ and $ l_{2,\mathfrak{z}}$ are determined as
\begin{equation}
	\begin{cases}
		l_{1,\mathfrak{z}} = w_{\mathfrak{z}}, \quad l_{2,\mathfrak{z}} = h_{\mathfrak{z}} \quad \text{if } \Delta_\mathrm{o} = 1 - \delta_{\mathrm{o},\mathfrak{d},\mathfrak{z}} \\
		l_{1,\mathfrak{z}} = h_{\mathfrak{z}}, \quad l_{2,\mathfrak{z}} = w_{\mathfrak{z}} \quad \text{if } \Delta_\mathrm{o} = \delta_{\mathrm{o},\mathfrak{d},\mathfrak{z}}.
	\end{cases}
\end{equation}  
These constraints, together with \eqref{eq: 2D-SAT-rec}, ensure that the clusters are tightly joined to form a continuous shape.

\subsection{Linearization Techniques}
This section provides the mixed-integer linearization techniques used in this paper.

\subsubsection{Binary AND operation}
For three binary variables $\delta_a,\delta_b,\delta_c$ with $\delta_c=\delta_a\cdot\delta_b$, we can define the binary AND operation through
\begin{equation}
	\begin{aligned}
		\delta_c &\leq \delta_a, \\
		\delta_c &\leq \delta_b, \\
		\delta_c &\geq \delta_a + \delta_b - 1.
	\end{aligned}
	\label{eq: 2D-relax-b-AND}
\end{equation}   

\subsubsection{Trigonometric expression}
Given an angle $\theta$, the trigonometric values of the angle can be expressed as a linear combination of discrete values. For example, the cosine of the angle can be calculated as
\begin{equation}
	\cos(\theta) = \sum_{k=0}^{N_a - 1} \delta_{k} \cdot \cos(\theta_k),  
\end{equation}
where $\delta_k$ is the binary angle decision variable, and $\theta_k$ is the discretized angle. 
For trigonometric products like $\cos(\theta_1)\cdot\cos(\theta_2)$, we expand the formulation and relax it using the previous binary AND.

\subsubsection{Continuous variable times trigonometric value}
If a continuous variable $x$ is multiplied by $\cos(\theta)$, the expression can be linearized using the binary variables $\delta_k$ and a big-M formulation as follows:
\begin{equation}
	\begin{aligned}
		x_{\mathrm{c},k} &\leq x \cdot \cos(\theta_k) + M \cdot (1 - \delta_k), \\
		x_{\mathrm{c},k} &\geq x \cdot \cos(\theta_k) - M \cdot (1 - \delta_k), \\
		x_{\mathrm{c},k} &\leq M \cdot \delta_k, \\
		x_{\mathrm{c},k} &\geq -M \cdot \delta_k,
	\end{aligned}
\end{equation}    
where $x_{\mathrm{c},k}$ is a continuous variable. This allows us to linearize the multiplication of continuous variable with trigonometric value as
\begin{equation}
	x \cdot \cos(\theta) = \sum_{k=0}^{N_a - 1} x_{\mathrm{c},k}.
\end{equation}
We implement this technique to linearize the vectors multiplication in~\eqref{eq: 2D-SAT-rec} -- \eqref{eq: 2D-SAT-circle}, and \eqref{eq: 2D-relative-placement}.

\subsubsection{Absolute Disjunction Function} 
For an absolute disjunction function $|a| \geq b, \quad b \geq 0$, we can introduce a linearization using two binary variables $\delta_\mathrm{sgn}$ and $\delta_\mathrm{suc}$, along with a big-M formulation. The linearized expression becomes
\begin{equation}
	\begin{aligned}
		a &\geq b - M \cdot (1 - \delta_\mathrm{sgn}) - M \cdot (1 - \delta_\mathrm{suc}), \\
		a &\leq -b + M \cdot \delta_\mathrm{sgn} + M \cdot (1 - \delta_\mathrm{suc}), \\
		a &\leq b + M \cdot \delta_\mathrm{suc}, \\
		a &\geq -b - M \cdot \delta_\mathrm{suc}.
	\end{aligned}
\end{equation}
In this formulation, the decision variable $\delta_\mathrm{suc}$ indicates whether the condition is satisfied, and the sign variable $\delta_\mathrm{sgn}$ governs the sign of $a$. We leverage this technique to relax SAT formulation~\eqref{eq: 2D-SAT-rec} and \eqref{eq: 2D-SAT-circle}. For example, we relax~\eqref{eq: 2D-SAT-rec} by introducing four decision variables, subject to the constraint $\sum_{i=1}^{4}\delta_\mathrm{suc_i} \geq 1$,
ensuring at least one SAT condition is met. 

\subsection{Objective and Problem Formulation}
\begin{figure}[t]
	\centering
	\includegraphics[trim=0cm 0.5cm 0cm 1.0cm,clip, width=0.8\columnwidth]{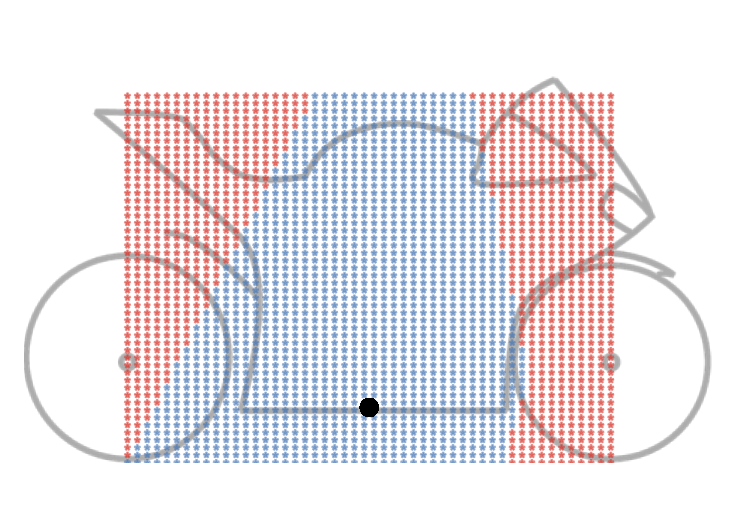}
	\caption{The CoG inactive region of the rear-wheel driven single \gls{acr:MM} topology, with the ideal \gls{acr:CoG} shown in black. The blue region represents the area where the maximum tractive force boundaries are inactive constraints, while the red region indicates where they are active constraints.}
	\label{fig: cog-dis}
\end{figure}

At the system level, the spatial arrangement of components directly determines the overall \gls{acr:CoG}, which in turn is a key factor governing the vehicle's handling. 
In extreme cases, the system's \gls{acr:CoG} can also influence the energy consumption. For instance, the tractive forces generated by the tires should be sufficient to accompany the requested torque from the rider. Yet the maximum tractive forces are limited by the friction coefficient and tire load through
\begin{equation}
	\begin{aligned}
		&F_\mathrm{t,max,f}(t) = \mu_\mathrm{f} \cdot F_\mathrm{n,f}(t),\\
		&F_\mathrm{t,max,r}(t) = \mu_\mathrm{r} \cdot F_\mathrm{n,r}(t),\\
	\end{aligned}
	\label{eq: ec-ft}
\end{equation}
where $F_{\mathrm{t,max},i}(t)$ with $i\in \{\text{f,r}\}$ represents the maximum tractive force for the front and rear tire, respectively, $\mu_i$ are the tire friction coefficients, and $F_{\mathrm{n},i}(t)$ are the normal forces acting on the tires. By disregarding suspension effects and taking the moment balance around each of the wheel centers, we can formulate the normal forces as
\begin{equation}
	\begin{aligned}
		& F_\mathrm{n,f}(t) = \frac{b}{l_\mathrm{wb}} \cdot m \cdot g - \frac{h}{l_\mathrm{wb}} \cdot m \cdot \dot{v}(t), \\ 
		& F_\mathrm{n,r}(t) = \frac{l_\mathrm{wb}-b}{l_\mathrm{wb}} \cdot m \cdot g + \frac{h}{l_\mathrm{wb}} \cdot m \cdot \dot{v}(t), \\ 
	\end{aligned}
	\label{eq: ec-n_force}
\end{equation}
where $b$ is the horizontal length from the \gls{acr:CoG} to the rear wheel, $h$ is the vertical length from the \gls{acr:CoG} to the ground, $l_\mathrm{wb}$ is the wheelbase length, $m$ is the total mass of the vehicle including the rider, and $\dot{v}(t)$ is the longitudinal acceleration. 
The total mass $m$ is defined as
\begin{equation}
	m = m_\mathrm{c} + m_\mathrm{r} + \sum_{\tau_j \in \mathbf{T}} {m_{\tau_j}},
	\label{eq: ec-sum-m}
\end{equation}
where $m_\mathrm{c}$ is the chassis mass, $m_\mathrm{r}$ is the rider mass, and $m_{\tau_j}$ represents the mass of each component that is present in topology $\mathbf{T}$.
Considering that the \gls{acr:CoG} is defined by the point $\left(b,h\right)$, it has a direct influence on the grip limits. Hence, we require the \gls{acr:CoG} to be positioned such that the normal forces are not limiting the tractive torque under normal riding conditions, both in positive and negative direction. If the objective is to minimize energy consumption, the ability to recuperate energy during deceleration is especially important. To check which \gls{acr:CoG} locations result in the maximum tractive forces being active constraints, we evaluate the tractive forces over the World Motorcycle Test Cycle (WMTC), which is part of the energy consumption minimization step in the integrated framework shown in Fig.~\ref{fig: Fig1_demo}. This evaluation for a grid of \gls{acr:CoG} locations on a rear-wheel driven motorcycle with single \gls{acr:MM} configuration is shown in Fig.~\ref{fig: cog-dis}. For dual-motor configurations, we compute the power-split between the motors based on efficiency maps, such that the overall motor losses are minimized. To maximize handling performance, we choose a single point $(x_{\text{ideal}},y_{\text{ideal}})$ within the inactive region as the ideal \gls{acr:CoG}, which is in this case the lowest point in the middle of the wheels.
The overall center of gravity $(x_{\text{CoG}},y_{\text{CoG}})$ is calculated as
\begin{equation}
	\begin{bmatrix} x_{\text{CoG}}\\ y_{\text{CoG}}\end{bmatrix}
	=
	\frac{\left(\sum_{\tau_j \in\mathbf{T}} \begin{bmatrix} x_{\tau_j}\\ y_{\tau_j}\end{bmatrix} m_{\tau_j}\right) + \begin{bmatrix} x_\mathrm{c}\\ y_\mathrm{c}\end{bmatrix} m_\mathrm{c} + \begin{bmatrix} x_r\mathrm{}\\ y_\mathrm{r}\end{bmatrix} m_\mathrm{r}}
	{\sum_{\tau_j\in\mathbf{T}} m_{\tau_j} + m_\mathrm{c} + m_\mathrm{r}}.
	\label{eq:cog}
\end{equation}


Finally, we define the \gls{acr:CoG} positioning objective
\begin{equation}
	J_{\text{dis}} = \frac{1}{l_n} \begin{Vmatrix} x_{\text{CoG}} - x_{\text{ideal}}\\ y_{\text{CoG}} - y_{\text{ideal}}\end{Vmatrix}^2_2,
\end{equation}
where $l_n$ is a normalization length. 
In case of a \gls{acr:MM} configuration, we include an additional constraints and objective for the \gls{acr:MM} placement and chain length:
\begin{equation}
	\begin{aligned}
		& x_\mathrm{MM} \in [\underline{x}_\mathrm{MM}, \overline{x}_\mathrm{MM}] 
		, \quad y_\mathrm{MM} \in [\underline{y}_\mathrm{MM}, \overline{y}_\mathrm{MM}], \\
		&J_\mathrm{dis,MM} = \frac{1}{l_\mathrm{n,MM}} \cdot 
		\begin{Vmatrix}
			x_\mathrm{MM}-x_\mathrm{w,r} \\ 
			y_\mathrm{MM}-y_\mathrm{w,r}
	\end{Vmatrix}^2_2,
\end{aligned}
	\label{eq: 2D-obj-ndis-MM}
\end{equation}    
where $[\underline{x}_\mathrm{MM}, \overline{x}_\mathrm{MM}]$ and $ [\underline{y}_\mathrm{MM}, \overline{y}_\mathrm{MM}]$ are the boundaries for the $\mathrm{MM}$, $J_\mathrm{dis,MM}$ is the non-negative normalized distance from $\mathrm{MM}$ to the rear wheel $\mathrm{w,r}$, and $l_\mathrm{n,MM}$ is the normalization length for $\mathrm{MM}$. 

By defining the sum of these two normalized distances as the handling performance objective, where $J_\mathrm{dis,MM}$ is weighted by the $\mathrm{MM}$ existence variable $\mathbf{n}_\mathrm{MM}$, we frame the problem as follows:
\begin{prob}[Optimal component placement]\label{prob: 2D placement}
	The optimal 2D component placement is the solution of
	\begin{equation*}
		\begin{aligned}
			&\min \: J_\mathrm{dis} + J_\mathrm{dis,MM}\cdot \mathbf{n}_\mathrm{MM}, \\
			\textnormal{s.t. } \:  &\eqref{eq: 2D-d-angle}-\eqref{eq: 2D-obj-ndis-MM}.
		\end{aligned}
	\end{equation*}   
\end{prob}
\noindent
Problem 1 can be solved with mixed-integer quadratic programming solvers with global optimality guarantees upon convergence~\cite{Lee2012}.

%
%
%
%
%

\section{Results}
\label{sec:results}
This section demonstrates our placement framework for two different electric motorcycle topologies. 
The first case includes a conventional single \gls{acr:MM} configuration, where the \gls{acr:em} is located inside the frame and connected to the rear wheel through a chain, as illustrated in Fig.~\ref{fig: MM-2D-tp}. The battery pack is modeled as a subsystem containing internal modules and we want to jointly optimize its geometry together with the placement of the other elements. Considering the application where space is typically very limited, accounting for freedom in the subsystem geometry design is a crucial aspect. Yet solving the full problem, whereby each internal module is a cluster, would result in impractical computation times. Therefore, we maintain tractable computation times by gradually increasing battery pack design complexity $N_{\text{com,BP}}$ and investigate how this complexity affects element placement.
The second case features a dual motor configuration consisting of a \gls{acr:MM} and a rear \gls{acr:HM}, where the latter is located inside the rear wheel. Similarly to the previous case, the battery pack is a subsystem and we compute the optimal component placement for various complexity levels.   

The mixed-integer optimization problems are solved using Gurobi~\cite{GurobiOptimization2021} and performed on an AMD Ryzen 7 5800H processor 3.20 GHz with 16 GB of RAM. Thereby, the computation time was approximately \unit[10]{min} to reach \unit[0\%]{optimality gap}. 



\subsection{Single motor configuration}
\begin{figure}[t] 
	\centering
	\begin{subfigure}[b]{0.49\linewidth} 
		\includegraphics[trim= 0.5cm 0.5cm 0.5cm 0.5cm, clip,width=\linewidth]{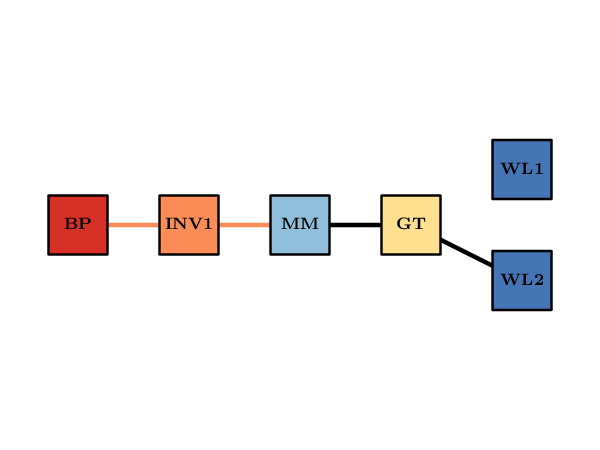}
		\caption{Single \gls{acr:MM} configuration}
		\label{fig: MM-2D-tp}
	\end{subfigure}  
	\hfill 
	\begin{subfigure}[b]{0.49\linewidth}
		\includegraphics[trim= 0.5cm 0.5cm 0.5cm 0.5cm, clip,width=\linewidth]{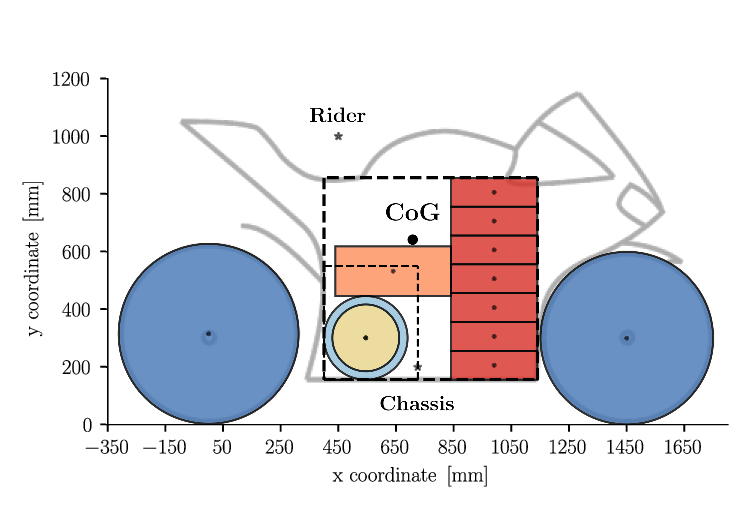}
		\caption{Benchmark solution}
		\label{fig: MM-2D-1com}
	\end{subfigure}
	
	\begin{subfigure}[b]{0.49\linewidth}
		\includegraphics[trim= 0.5cm 0.5cm 0.5cm 0cm, clip,width=\linewidth]{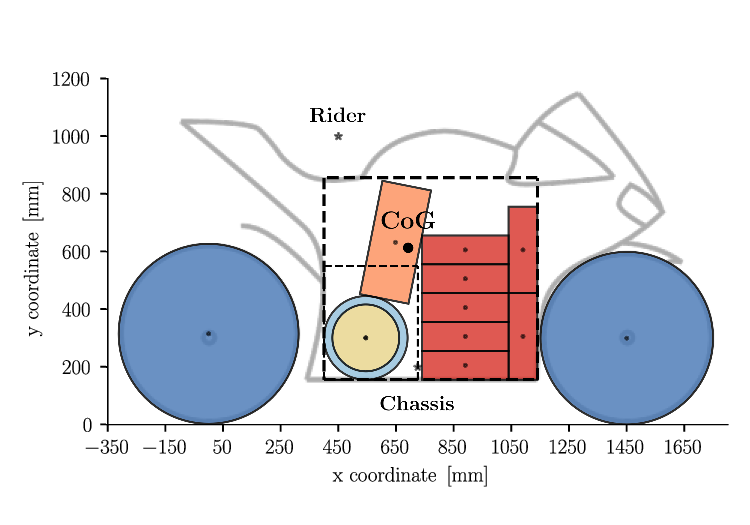}
		\caption{$N_{\text{com,BP}}=2$ solution}
		\label{fig: MM-2D-2com}
	\end{subfigure}
	\hfill
	\begin{subfigure}[b]{0.49\linewidth}
		\includegraphics[trim= 0.5cm 0.5cm 0.5cm 0cm, clip,width=\linewidth]{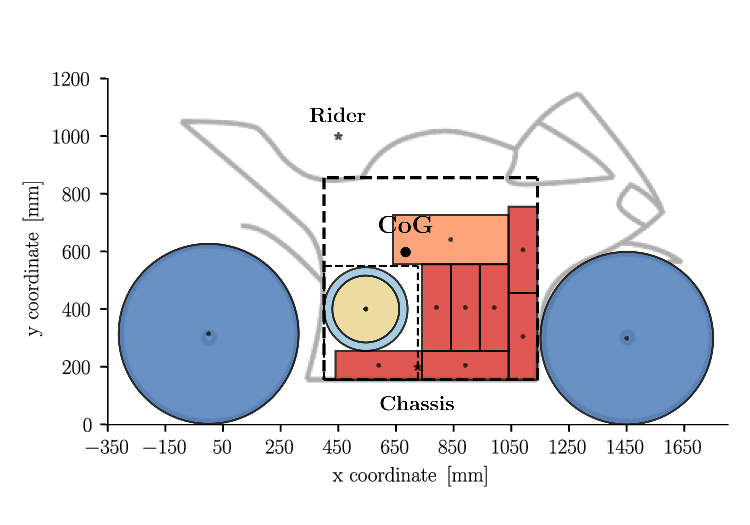}
		\caption{$N_{\text{com,BP}}=3$ solution}
		\label{fig: MM-2D-3com}
	\end{subfigure}
	\caption{A single $\mathrm{MM}$ topology and its optimal 2D placements with different complexity. The outer dashed black line represents the design space, while the inner dashed black line encloses the feasible region for the mounted motor. The chassis and rider CoGs are marked with $\star$, the overall CoG is represented by a solid black circle, and the CoGs of components and sub-modules are indicated by black dots.}
	\label{fig: MM-2D}
\end{figure}

For the single \gls{acr:MM} configuration, we consider three cases with varying complexity in placement of the modules inside the battery pack. 
As a benchmark, we consider a battery pack complexity of $N_{\text{com,BP}}=1$, which is shown in Fig.~\ref{fig: MM-2D-1com}. Here, all modules are forced to be placed in a single rectangular cluster. In this case, the optimal design of the battery pack is a single column stacked configuration. 
If we increase the complexity to $N_{\text{com,BP}}=2$, we obtain a battery pack configuration with two columns, as shown in Fig.~\ref{fig: MM-2D-2com}. In addition, we observe a non-orthogonal orientation of the inverter, demonstrating the importance of modeling finer orientations when space is limited. 
Lastly, we set the complexity $N_{\text{com,BP}}=3$ to obtain a battery pack with three clusters (two $1\times 2$ and one $1\times 3$). By increasing the complexity, we obtain the overall best placement, with a highly irregular geometry of the battery pack and an overall improvement of \unit[2.5]{\%} in the objective, compared to the benchmark solution. This highlights the importance of irregular geometries in confined spaces.


\subsection{Dual motor configuration}

\begin{figure}[t]
	\begin{center}
		\begin{subfigure}[b]{0.25\linewidth}
			\includegraphics[trim= 1.5cm 0.5cm 1.5cm 0.5cm, clip,width=\linewidth]{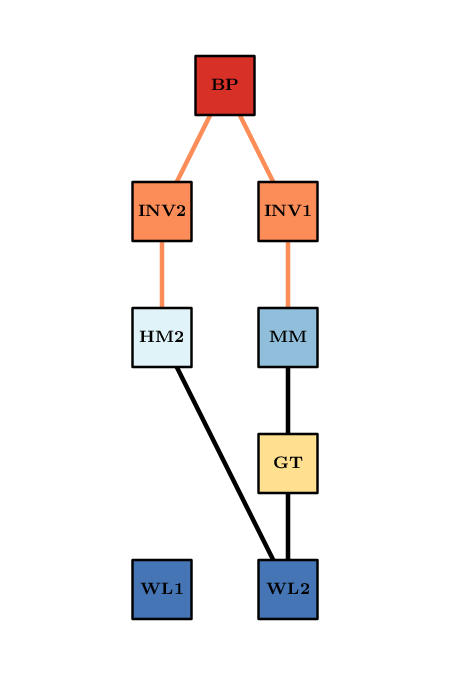}
		\end{subfigure}
		\hfill 
		\begin{subfigure}[b]{0.73\linewidth}
			\includegraphics[trim= 0.5cm 0.5cm 0.5cm 0.5cm, clip,width=\linewidth]{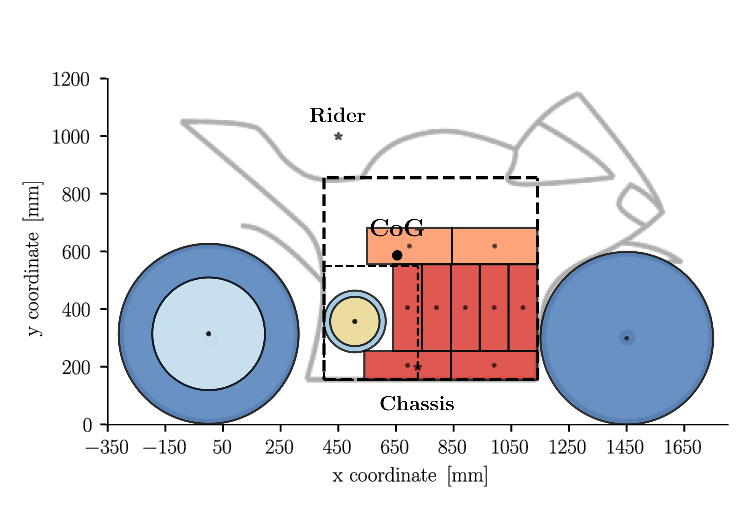}
		\end{subfigure}
		\caption{The topology with a $\mathrm{MM}$ and a $\mathrm{HM_2}$, together with the placement results for battery pack complexity $N_\mathrm{com,BP}$ of 2 and 3.} 
		\label{fig: HM2-MM-2D}
	\end{center}
\end{figure}

In the second configuration, we introduce a \gls{acr:HM} and an additional inverter, as shown in the topology in Fig.~\ref{fig: HM2-MM-2D}. To retain the same system power output, we scale the \gls{acr:MM} and inverters by a factor of 0.5, compared to the single motor case. In this case, battery pack complexities $N_{\text{com,BP}}=2$ and $3$ yield the same optimal placement, indicating that additional design freedom does not always lead to improved solutions. This confirms that we can simplify the placement problem by gradually increasing subsystem complexity until no further improvements are observed, ultimately leading to more manageable computation times.   

\section{Conclusion} \label{sec: conclusion}
In this paper, we proposed a two-dimensional optimization framework that determines the powertrain components’ placement in electric motorcycles, while maintaining optimal energy consumption and handling performance. To address the main challenge of limited packaging space in motorcycles, our method introduces the flexibility to manage non-orthogonal orientations and both rectangular and circular shapes, while guaranteeing mixed-integer convexity. 
Additionally, we devised a method that enables irregular geometries of subsystems through the introduction of clusters within each subsystem. For the studied use cases, this has shown to be beneficial in bringing the overall \gls{acr:CoG} \unit[2.5]{\%} closer to the ideal \gls{acr:CoG} and has the additional benefit of maintaining the computation time at a manageable level.  
Future work may extend this framework to additional vehicle topologies, refine the optimization process for larger systems, or apply it to other engineering applications where component placement is critical.

\section{Acknowledgment}
\noindent We thank Dr.~I.~New for proofreading this paper. This paper was partly supported by the NEON research project (project number 17628 of the Crossover program which is (partly) financed by the Dutch Research Council (NWO)).


\bibliographystyle{IEEEtran}
\bibliography{main,SML_papers,powertrains}

\newcommand{\noopsort}[1]{} \newcommand{\printfirst}[2]{#1}
  \newcommand{\singleletter}[1]{#1} \newcommand{\switchargs}[2]{#2#1}
\begin{thebibliography}{10}
\providecommand{\url}[1]{#1}
\csname url@samestyle\endcsname
\providecommand{\newblock}{\relax}
\providecommand{\bibinfo}[2]{#2}
\providecommand{\BIBentrySTDinterwordspacing}{\spaceskip=0pt\relax}
\providecommand{\BIBentryALTinterwordstretchfactor}{4}
\providecommand{\BIBentryALTinterwordspacing}{\spaceskip=\fontdimen2\font plus
\BIBentryALTinterwordstretchfactor\fontdimen3\font minus
  \fontdimen4\font\relax}
\providecommand{\BIBforeignlanguage}[2]{{%
\expandafter\ifx\csname l@#1\endcsname\relax
\typeout{** WARNING: IEEEtran.bst: No hyphenation pattern has been}%
\typeout{** loaded for the language `#1'. Using the pattern for}%
\typeout{** the default language instead.}%
\else
\language=\csname l@#1\endcsname
\fi
#2}}
\providecommand{\BIBdecl}{\relax}
\BIBdecl

\bibitem{WijknietHofman2018}
J.~Wijkniet and T.~Hofman, ``Modified computational design synthesis using
  simulation-based evaluation and constraint consistency for vehicle powertrain
  systems,'' \emph{{IEEE Trans.\ on Vehicular Tech.\ }}, vol.~67, no.~9, pp.
  8065--8076, 2018.

\bibitem{ChenouardHartmannEtAl2016}
R.~Chenouard, C.~Hartmann, A.~Bernard, and E.~Mermoz, ``Computational design
  synthesis using model-driven engineering and constraint programming,''
  \emph{{Lecture Notes in Computer Science}}, 2016.

\bibitem{ZhaoTangEtAl2022}
Z.~Zhao, P.~Tang, and H.~Li, ``Generation, screening, and optimization of
  powertrain configurations for power-split hybrid electric vehicle: A
  comprehensive overview,'' \emph{{IEEE Transactions on Transportation
  Electrification}}, 2022.

\bibitem{KabalanVinotEtAl2020}
B.~Kabalan, E.~Vinot, R.~Trigui, and C.~Dumand, ``Systematic methodology for
  architecture generation and design optimization of hybrid powertrains,''
  \emph{{IEEE Trans.\ on Vehicular Tech.\ }}, vol.~69, no.~12, pp.
  14\,846--14\,857, 2020.

\bibitem{BerxGadeyneEtAl2014}
K.~Berx, K.~Gadeyne, M.~Dhadamus, G.~Pipeleers, and G.~Pinte, ``Model-based
  gearbox synthesis,'' in \emph{{Mechatronics Forum International Conference}},
  2014.

\bibitem{PiacentiniCheongEtAl2020}
C.~Piacentini, H.~Cheong, M.~Ebrahimi, and A.~Butscher, ``Multi-speed gearbox
  synthesis using global search and non-convex optimization,'' in
  \emph{{Integration of Constraint Programming, Artificial Intelligence, and
  Operations Research}}, 2020.

\bibitem{KampenSalazarEtAl2023}
J.~van Kampen, M.~Salazar, and T.~Hofman, ``A two-dimensional spatial
  optimization framework for vehicle powertrain systems,'' in \emph{{IEEE
  Vehicle Power and Propulsion Conference}}, 2023.

\bibitem{KampenSalazarEtAl2025}
------, ``Automated three-dimensional spatial optimization for multi-domain
  systems with alignment constraints,'' \emph{{ASME Journal of Mechanical
  Design}}, 2025, in Press.

\bibitem{HofstetterLechleitnerEtAl2018}
M.~Hofstetter, D.~Lechleitner, M.~Hirz, M.~Gintzel, and A.~Schmidhofer,
  ``{Multi-objective gearbox design optimization for xEV-axle drives under
  consideration of package restrictions},'' \emph{{Forschung im
  Ingenieurwesen}}, vol.~82, pp. 361--370, 2018.

\bibitem{FilomenoAhmadEtAl2021}
G.~Filomeno, M.~A. Ahmad, M.~Wolkerstorfer, B.~Krueger, P.~Tenberge, and
  D.~Dennin, ``Multi-objective electric powertrain design optimization under
  package constraints,'' in \emph{{Electromechanical Drive Systems}}, 2021.

\bibitem{Laignel2021}
G.~Laignel, N.~Pozin, X.~Geffrier, L.~Delevaux, F.~Brun, and B.~Dolla, ``Floor
  plan generation through a mixed constraint programming-geneticoptimization
  approach,'' \emph{{Automation in Construction}}, 2021.

\bibitem{Tresca2023}
G.~Tresca, G.~Cavone, R.~Carli, and M.~Dotoli, ``A mathematical model for the
  optimal configuration of automatedstorage systems with sliding trays,'' in
  \emph{{European Control Conference}}, 2023.

\bibitem{CaganShimadaEtAl2002}
J.~Cagan, K.~Shimada, and S.~Yin, ``A survey of computational approaches to
  three-dimensional layout problems,'' \emph{{Computer-aided design}}, 2002.

\bibitem{BelloEtAl2024}
W.~B. Bello, S.~R.~T. Peddada, A.~Bhattacharyya, L.~E. Zeidner, J.~T. Allisona,
  and K.~A. James, ``Multi-physics three-dimensional component placement and
  routing optimization using geometric projection,'' \emph{{ASME Journal of
  Mechanical Design}}, 2024.

\bibitem{PeddadaEtAl2022}
S.~R.~T. Peddada, L.~E. Zeidner, H.~T. Ilies, K.~A. James, and J.~T. Allison,
  ``Toward holistic design of spatial packaging of interconnected systems with
  physical interactions (spi2),'' \emph{{ASME Journal of Mechanical Design}},
  2022.

\bibitem{Liang2015}
C.~Liang and X.~Liu, ``The research of collision detection algorithm based on
  separating axis theorem,'' \emph{{Int.\ Journal of Science}}, 2015.

\bibitem{RichardsHow2005}
A.~Richards and J.~How, ``Mixed-integer programming for control,'' in
  \emph{{Proc.\ of the American Control Conference}}, 2005.

\bibitem{Lee2012}
J.~Lee and S.~Leyffer, Eds., \emph{Mixed Integer Nonlinear Programming}.\hskip
  1em plus 0.5em minus 0.4em\relax {Springer-Verlag}, 2012.

\bibitem{GurobiOptimization2021}
{Gurobi Optimization, LLC}. (2021) Gurobi optimizer reference manual. Available
  online at \url{http://www.gurobi.com}.

\end{thebibliography}


\end{document}